

  \newcount\fontset


                   \fontset=2


  \def\dualfont#1#2#3{\font#1=\ifnum\fontset=1 #2\else#3\fi}
  \dualfont\eightrm {cmr8} {cmr7}
  \dualfont\eightsl {cmsl8} {cmr7}
  \dualfont\eightit {cmti8} {cmti10}
  \dualfont\eightmi {cmmi8} {cmmi7}
  \dualfont\tensc {cmcsc10} {cmcsc10}
  \dualfont\titlefont {cmbx12} {cmbx10}
  \dualfont\eightsymbol {cmsy8} {cmsy10}


  \magnification=\magstep1
  \nopagenumbers
  \voffset=2\baselineskip
  \advance\vsize by -\voffset
  \headline{\ifnum\pageno=1 \hfil \else \tensc\hfil
    group and inverse semigroup actions
  \hfil\folio \fi}


  \def\vg#1{\ifx#1\null\null\else
    \ifx#1\ { }\else
    \ifx#1,,\else
    \ifx#1..\else
    \ifx#1;;\else
    \ifx#1::\else
    \ifx#1''\else
    \ifx#1--\else
    \ifx#1))\else
    { }#1\fi \fi \fi \fi \fi \fi \fi \fi \fi}
  \def\.{.\enspace}
  \newcount\secno \secno=0
  \newcount\stno
  \def\goodbreak{\vskip0pt plus.1\vsize\penalty-250
    \vskip0pt plus-.1\vsize\bigskip}
  \outer\def\section#1{\stno=0
    \global\advance\secno by 1
    \goodbreak\vskip\parskip
    \message{\number\secno.#1}
    \noindent{\bf\number\secno\.#1\.}}
  \def\state#1 #2\par{\advance\stno by 1\medbreak\noindent
    {\bf\number\secno.\number\stno.\enspace #1.\enspace}{#2}\medbreak}
  \def\nstate#1 #2#3\par{\state{#1} {#3}\par
    \edef#2{\number\secno.\number\stno}}
  \def\proof{\medbreak\noindent{\it Proof\.}}
  \def\proofend{\ifmmode\eqno\square\else\hfill\square\medbreak\fi}
  \newcount\zitemno \zitemno=0
  \def\zitem{\global\advance\zitemno by 1 \smallskip
    \item{\ifcase\zitemno\or i\or ii\or iii\or iv\or v\or vi\or vii\or
    viii\or ix\or x\or xi\or xii\fi)}}
  \def\$#1{#1 $$$$ #1}


  \newcount\bibno \bibno=0
  \def\newbib#1{\advance\bibno by 1 \edef#1{\number\bibno}}
  \def\cite#1{{\rm[\bf #1\rm]}}
  \def\scite#1#2{{\rm[\bf #1\rm, #2]}}
  \def\lcite#1{(#1)}
  \def\se#1#2#3#4{\def\a{#1}\def\b{#2}\ifx\a\b#3\else#4\fi}
  \def\setem#1{\se{#1}{}{}{, #1}}
  \def\index#1{\smallskip \item{[#1]}}

  \def\ATarticle#1{\zarticle #1 xyzzy }
  \def\zarticle#1, author = #2,
   title = #3,
   journal = #4,
   year = #5,
   volume = #6,
   pages = #7,
   NULL#8 xyzzy {\index{#1} #2, ``#3'', {\sl #4\/} {\bf #6} (#5), #7.}

  \def\ATtechreport#1{\ztechreport #1 xyzzy }
  \def\ztechreport#1,
    author = #2,
    title = #3,
    institution = #4,
    year = #5,
    note = #6,
    toappear = #7,
    NULL#8
    xyzzy
    {\index{#1} #2, ``#3''\setem{#6}\setem{#4}\setem{#5}\se{#7}{}{.}{,
to appear in #7}}

  \def\ATbook#1{\zbook #1 xyzzy }
  \def\zbook#1,
    author = #2,
    title = #3,
    publisher = #4,
    year = #5,
    volume = #6,
    series = #7,
    NULL#8
    xyzzy
    {\index{#1} #2, ``#3''\setem{#7}\se{#6}{}{}{ vol. #6}, #4, #5.}

  \def\ATmasterthesis#1{\zmasterthesis #1 xyzzy }
  \def\zmasterthesis#1,
    author = #2,
    title = #3,
    school = #4,
    year = #5,
    NULL#6
    xyzzy
    {\index{#1} #2, ``#3'', Master Thesis, #4, #5.}


  \def\stress#1{{\it #1}\/}
  
  \def\crossproduct{\hbox to 1.8ex{$\times \kern-.45ex\vrule height1.1ex
depth0pt width0.45truept$\hfill}}
  \def\:{\colon}
  \def\*{\otimes}
  \def\+{\oplus}
  \def\x{\times}
  \def\({\bigl(}
  \def\){\bigl)}
  \def\arw{\rightarrow}
  \def\cstar{$C^*$}
  \def\square{\hbox{$\sqcap\!\!\!\!\sqcup$}}
  \def\for#1{,\quad #1}


  \def\sg{semi\-group\vg}
  \def\sgs{semi\-groups\vg}
  \def\isg{inverse \sg}
  \def\isgs{inverse \sgs}
  \def\homo{ho\-mo\-mor\-phism\vg}
  \def\homos{ho\-mo\-mor\-phisms\vg}
  \def\I{{\cal I}}
  \def\S#1{{\cal S}(#1)}
  \def\dom{\hbox{\rm dom}}
  \def\ran{\hbox{\rm ran}}
  \def\directsum{\bigoplus}
  \def\t{\theta}
  \def\e{\varepsilon}
  \def\p{P}
  \def\q{Q}
  \def\span{\overline{\hbox{\rm span}}}
  \def\degree{\partial}
  \def\inv{^{-1}}
  \def\Pe{{\cal P}_e(G)}
  \def\C{{\bf C}}
  \def\B{{\cal B}(H)}
  
  \def\Pl{{\cal P}_\ell(A)}
  \def\F{{\cal F}\bigl({\cal P}_e(G)\bigr)}


  \newbib{\AEE}
  \newbib{\Barnes}
  \newbib{\DP}
  \newbib{\ExelTPA}
  \newbib{\ExelCircle}
  \newbib{\FD}
  \newbib{\Howie}
  \newbib{\Jensen}
  \newbib{\McClanahan}
  \newbib{\Pedersen}
  \newbib{\Pimsner}
  \newbib{\QuiggRaeburn}
  \newbib{\Renault}
  \newbib{\Sieben}


  \null
  \vskip-2\bigskipamount

  \begingroup
  \def\c{\centerline}
  \c{\titlefont PARTIAL ACTIONS OF GROUPS AND ACTIONS}
  \medskip
  \c{\titlefont OF INVERSE SEMIGROUPS}

  \bigskip
  \baselineskip=10pt
  \eightit
  \c{\tensc Ruy Exel\footnote{*}{\eightrm Partially supported by CNPq,
Brazil.}}
  \c{Departamento de Matem\'atica}
  \c{Universidade de S\~ao Paulo}
  \c{Rua do Mat\~ao, 1010}
  \c{05508-900 S\~ao Paulo -- Brazil}
  \c{exel@ime.usp.br}

  %
  \def\G{{\eightmi G}\vg}
  
  \def\SG{{\eightsymbol S}({\eightmi G})\vg}
  \bigskip\bigskip\eightrm\baselineskip=3ex
  \midinsert\narrower\narrower
  Given a group \G, we construct, in a canonical way, an \isg\ \SG
associated to \G.  The actions of \SG are shown to be in one-to-one
correspondence with the partial actions of \G, both in the case of
actions on a set, and that of actions as operators on a Hilbert space.
In other words, \G and \SG have the same representation theory.

We show that \SG governs the sub\sg of all closed linear subspaces of a
\G-graded {\eightmi C}*-algebra, generated by the grading subspaces.  In
the special case of finite groups, the maximum number of such subspaces
is computed.

A ``partial'' version of the group {\eightmi C}*-algebra of a discrete
group is introduced.  While the usual group {\eightmi C}*-algebra of
finite commutative groups forgets everything but the order of the group,
we show that the partial group {\eightmi C}*-algebra of the two
commutative groups of order four, namely {\eightmi Z}/4{\eightmi Z} and
{\eightmi Z}/2{\eightmi Z} $\oplus$ {\eightmi Z}/2{\eightmi Z}, are not
isomorphic.
  \endinsert \endgroup

  \section {Introduction} A \sg\ $S$ is said to be an \isg\ provided
there exists, for each $x$ in $S$, a unique element $x^*$ in $S$ such
that
  \zitemno=0
  \zitem $xx^*x = x$,
  \zitem $x^*xx^* = x^*$.

  \smallskip According to one of the earliest results on this subject,
the Vagner--Preston representation Theorem \scite{\Howie}{V.1.10}, every
\isg is isomorphic to a sub\sg of $\I(X)$, the \isg of partially defined
bijective maps on a set $X$.  The interested reader will also find in
\cite{\Howie} the basic facts about \isgs.

By a partially defined map on $X$ we mean a map $\phi \colon A \arw B$,
where $A$ and $B$ are subsets of $X$.  Given such a $\phi$ we will
denote by $\dom(\phi)$ and $\ran(\phi)$ its domain $A$ and range $B$,
respectively.  The multiplication rule on $\I(X)$ is given by
composition of partial maps in the largest domain where it makes sense,
that is, if $\phi$ and $\psi$ are elements of $\I(X)$ then
  $\dom(\phi\psi) = \psi\inv(\ran(\psi)\cap\dom(\phi))$
  and
  $\ran(\phi\psi) = \phi(\ran(\psi)\cap\dom(\phi))$.

The representation theory of \isgs studies \homos (i.e \sg-\homos) $\pi
\colon S \arw \I(X)$ where $S$ is the \isg under scrutiny and $X$ is a
set.  Adopting the point of view of the theory of group actions, we make
the following:

  \nstate Definition \ActionOfISG An \stress{action} of the \isg\ $S$ on
the set $X$ is a \homo\ $\pi \colon S \arw \I(X)$.

In view of the recent interest in the concept of partial actions of
groups on C*-algebras
  (see
  \cite{\ExelCircle},
  \cite{\McClanahan},
  \cite{\ExelTPA},
  \cite{\QuiggRaeburn},
  \cite{\Sieben},
  \cite{\Pimsner},
  \cite{\AEE}),
  it is interesting to study the more general concept of partial actions
of groups on sets, as defined below.

  \nstate Definition \PartialAction Given a group $G$ and a set $X$, a
\stress{partial action} of $G$ on $X$ is a pair
    $$
  \Theta = \(\{D_t\}_{t\in G}, \{\t_t\}_{t\in G}\),
    $$
  where, for each $t$ in $G$, $D_t$ is a subset of $X$ and $\t_t\colon
D_{t\inv} \arw D_t$ is a bijective map, satisfying, for all $r$ and $s$
in $G$
  \zitemno=0
  \zitem $D_e = X$ and $\t_e = id_X$, the identity map on $X$ (here, as
always, $e$ denotes the identity element of $G$),
  \zitem $\t_r(D_{r\inv}\cap D_s) = D_r\cap D_{rs}$,
  \zitem $\t_r(\t_s(x)) = \t_{rs}(x)
    \for{x\in D_{s\inv} \cap D_{s\inv r\inv}}$.

  Thus, a partial action of $G$ on $X$ is also a map from $G$ into
$\I(X)$, just like the case for actions of \isgs.  However, the map
$t\in G\mapsto \t_t\in \I(X)$ does not satisfy $\t_r\t_s = \t_{rs}$ as
one might imagine, but rather we have $\t_r\t_s \subseteq \t_{rs}$, that
is, $\t_{rs}$ is an extension of $\t_r\t_s$.

It is the goal of this work to introduce, for each group $G$, an \isg,
which we will denote by $\S G$, such that the partial actions of $G$ on
a set $X$ are in one-to-one correspondence with the \homos from $\S G$
into $\I(X)$.

This problem was first considered in Sieben's \cite{\Sieben} master
thesis, where a partial answer was given in the case of actions on
C*-algebras.  The \isg constructed there, however, is not intrinsically
obtained from $G$, but it depends also on the partial action upon
consideration.

  \section {The \isg associated to a group} Throughout this section $G$
will denote a fixed group.  The crucial property of partial actions
which will enable us to construct our \isg is that, although $\t_{rs}$
is only an extension of $\t_r\t_s$ (notation as in the previous
section), we always have
  $\t_r\t_s\t_{s\inv} = \t_{rs}\t_{s\inv}$.
  More will be said about this property later.

  \nstate Definition \DefSG We let $\S G$ denote the universal \sg
defined via generators and relations as follows.  To each element $t$ in
$G$ we take a generator $[t]$ (from any fixed set having as many
elements as $G$).  For every pair of elements $t$, $s$ in $G$ we
consider the relations
  \zitemno=0
  \zitem $[s\inv][s][t] = [s\inv][st]$,
  \zitem $[s][t][t\inv] = [st][t\inv]$,
  \zitem $[s][e] = [s]$,
  \zitem $[e][s] = [s]$.

A rule of thumb to remember the first two relations is that $[s][t]$ can
be replaced by $[st]$ only if there is, either an $[s\inv]$ hanging out
on the left hand side, or a $[t\inv]$ on the right.

Note that, as a consequence of (i) and (iii), we have that
$[t][t\inv][t]=[t]$, which is the first hint that $\S G$ might be shown
to be an \isg.  This will, in fact, be accomplished below.  We should
also remark that axiom (iv) is a consequence of the previous ones, since
we have, for each $s$ in $G$, that
  $$[e][s] = [ss\inv][s] = [s][s\inv][s] = [s].$$
  Therefore (iv) could be removed from the definition without any
significant consequence.

By the universal property of \sg\null s defined via generators and
relations we have:

  \nstate Proposition \Universal Given a \sg\ $S$ and a map $f\colon G
\arw S$ satisfying
  \zitemno=0
  \zitem $f(s\inv) f(s) f(t) = f(s\inv) f(st)$,
  \zitem $f(s) f(t) f(t\inv) = f(st) f(t\inv)$,
  \zitem $f(s) f(e) = f(s)$,
  \smallskip\noindent there exists a unique \homo\ $\tilde f \colon \S G
\arw S$ such that $\tilde f([t]) = f(t)$.

  The proof is elementary but we would like to point out that, according
to the remarks above, the relation $f(e) f(s) = f(s)$, corresponding to
\lcite{\DefSG.iv}, need not be included as it is a consequence of (i) --
(iii).

  \state Proposition \sl There is an involutive anti-automorphism
$*\colon \S G \arw \S G$ such that, for every $t$ in $G$,
$[t]^*=[t\inv]$.

  \proof Let $\S G^{op}$ be the opposite \sg, that is $\S G^{op}$
coincides with $\S G$ in all respects, except that multiplication of the
elements $\alpha$ and $\beta$ in $\S G^{op}$ is defined by
  $$\alpha \bullet \beta := \beta\alpha$$
  where $\beta\alpha$, on the right hand side, corresponds to the usual
multiplication in $\S G$.
  Define $f\colon G \arw \S G$ by $f(t) = [t\inv]$.  It is easy to see
that, for $s$ and $t$ in $G$ we have
  \zitemno=0
  \zitem $f(s\inv)\bullet f(s)\bullet f(t) = f(s\inv)\bullet f(st)$,
  \zitem $f(s)\bullet f(t)\bullet f(t\inv) = f(st)\bullet f(t\inv)$,
  \zitem $f(s)\bullet f(e) = f(s)$.
  \smallskip \noindent So, by \lcite{\Universal}, $f$ extends to a
\homo\ $*\colon \S G \arw \S G^{op}$, which, when viewed as a map from
$\S G$ to itself, rather than into the opposite \sg, is an
anti-automorphism satisfying the desired properties.
  \proofend

  An element $\e$ in a \sg is said to be an idempotent if $\e^2=\e$.
Idempotents in $\S G$ ares studied next.

  \state Proposition \sl For every $t$ in $G$ let $\e_t = [t][t\inv]$.
Then each $\e_t$ is a self-adjoint idempotent in the sense that $\e_t^*
= \e_t$. In addition, $\e_t$ and $\e_s$ commute for every $t$ and $s$ in
$G$.

  \proof The only part of the statement requiring some care is the
commutativity of $\e_t$ and $\e_s$.  Using (i) and (ii) of
\lcite{\DefSG} repeatedly, as well as the fact that $[t][t\inv][t]=[t]$
we have,
    $$
  \e_t\e_s =
  [t][t\inv][s][s\inv] =
  [t][t\inv s][s\inv] =
  [t][t\inv s][s\inv t][t\inv s][s\inv]
    \$ =
  [s][s\inv t][t\inv s][s\inv] =
  [s][s\inv t][t\inv] =
  [s][s\inv][t][t\inv] =
  \e_s\e_t. \proofend
    $$

  \nstate Proposition \Decomposition Every element $\alpha$ in $\S G$
admits a decomposition
    $$
  \alpha = \e_{r_1}\ldots \e_{r_n} [s],
    $$
  where $n\geq 0$ and $r_1, \ldots, r_n, s$ are elements of $G$.  In
addition, one can assume that
  \zitemno = 0
  \zitem $r_i \neq r_j$ for $i\neq j$,
  \zitem $r_i \neq s$ and $r_i \neq e$ for all $i$.

  \proof Let $S$ be the subset of $\S G$ consisting of those $\alpha$
that do admit a decomposition as above.  Since $n=0$ is allowed, we see
that each $[s]$ belongs to $S$.  To prove the statement it is then
enough to verify that $S$ is a sub\sg of $\S G$, in view of the fact
that the $[s]$'s form a generating system for $\S G$.

Let $\alpha = \e_{r_1}\ldots \e_{r_n}[s]$.  It suffices to demonstrate
that $\alpha [t]$ belongs to $S$, since this will prove $S$ to be a
right ideal and hence a sub\sg.  Now, note that
    $$
  [s][t] =
  [s][s\inv][s][t] =
  [s][s\inv][st] =
  \e_s[st]
    $$
  So
    $$
  \alpha [t] =
  \e_{r_1}\ldots \e_{r_n}[s][t] =
  \e_{r_1}\ldots \e_{r_n}\e_s[st] \in S.
    $$

  With regard to the last part of the statement, since the idempotents
$\e_{r_j}$ commute with each other, it is easy to see how to eliminate
repetitions among the $\e_r's$, if any, as well as any occurrence of
$\e_e$, which is the identity element of $\S G$.

Also, if some $r_i = s$ then $\e_{r_i} = [s][s\inv]$ and, again by
commutativity of idempotents, we have
    $$
  \alpha =
  \e_{r_1}\ldots\widehat {\e_{r_i}} \ldots e_{r_n} [s][s\inv][s] =
  \e_{r_1}\ldots\widehat {\e_{r_i}} \ldots e_{r_n} [s].
    $$
  So $\e_{r_i}$ can be eliminated.  \proofend

  \state Definition If $\alpha$ is written as $\alpha = \e_{r_1}\ldots
\e_{r_n} [s]$, in such a way that conditions (i) and (ii) of
\lcite{\Decomposition} are verified, we say that $\alpha$ is in
\stress{standard} form.

  We now present some more evidence to support our earlier claim (to be
proved below) that $\S G$ is a \isg.

  \state Proposition \sl For each $\alpha$ in $\S G$ one has
$\alpha\alpha^*\alpha = \alpha$ and $\alpha^*\alpha\alpha^* = \alpha^*$.

  \proof Let $\alpha$ be represented as $\alpha = \e_{r_1}\ldots
\e_{r_n}[s]$.  We then have
    $$
  \alpha \alpha^* \alpha =
  \e_{r_1}\ldots \e_{r_n}[s]
  [s\inv] \e_{r_n}\ldots \e_{r_1}
  \e_{r_1}\ldots \e_{r_n}[s] =
  \e_{r_1}\ldots \e_{r_n}[s][s\inv][s] =
  \alpha.
    $$

The second statement follows from the first one and our knowledge that
$*$ is an anti-automorphism of $\S G$.
  \proofend

  The only task remaining, before we can assert that $\S G$ is in fact
an \isg, is to show that $\alpha^*$ is the unique element satisfying the
identities in the statement of our previous result.  However, uniqueness
results in algebraic structures originated from generators and relations
are difficult to establish, unless we can find a separating family of
representations. This is the purpose of our next section.

  \section {Representations of $\S G$} In this section we shall use the
term ``representation'' in a very loose way, meaning any \homo of $\S G$
into a \sg.  These representations will often be obtained with the aid
of \lcite{\Universal}.

The most obvious representation of $\S G$ one can think of, is the map
$$ \degree\colon \S G \arw G $$ given by $\degree([t]) = t$.  We shall
refer to $\degree(\alpha)$ as the \stress{degree} of $\alpha$.  Clearly
$\degree(\e_{r_1}\ldots \e_{r_n} [s]) = s$.

  We next discuss a more subtle representation of $\S G$.  Let $\Pe$ be
the set of all finite subsets of $G$ which contain the unit element $e$.
Thus $G$ and $\{e\}$ are, respectively, the biggest and the smallest
elements of $\Pe$.  The representation we are about to consider is a
\homo from $\S G$ into the \sg $\F$ consisting of all functions from
$\Pe$ into itself, under the composition rule.

Thus, for each $t\in G$ let us denote by $\phi_t$ the map
    $$
  \phi_t\colon \Pe \arw \Pe
    $$
  given by $\phi_t(E) = tE \cup \{e\}$, for each $E\in \Pe$.
  Here, as usual, $tE$ stands for the set $\{ts: s\in E\}$.

  Observe that, since $e\in E$, we can also write $\phi_t(E) = tE \cup
\{t,e\}$.

  \state Proposition \sl The map $t\in G \mapsto \phi_t \in \F$
satisfies the properties (i) -- (iii) of \lcite{\Universal} and hence
there is a unique representation $\Lambda \colon \S G \arw \F$ such that
$\Lambda ([t]) = \phi_t$.

  \proof We prove only \lcite{\Universal.i}.  For every $t,s\in G$ and
each $E\in \Pe$, we have
      $$
    \phi_{s\inv} \phi_s \phi_t (E)
    = \phi_{s\inv} \phi_s (tE \cup \{t,e\})
    = \phi_{s\inv} (stE \cup \{st,s,e\})
    = tE \cup \{t,e,s\inv\},
      $$
  while
      $$
  \phi_{s\inv} \phi_{st} (E) =
  \phi_{s\inv} (stE \cup \{st,e\}) =
  tE \cup \{t,s\inv,e\}. \proofend
      $$

  Observe that, for any $r$ in $G$ and $E$ in $\Pe$, one has
$\Lambda(\e_r)(E) = E \cup \{r\}$.  In particular, if $\alpha =
\e_{r_1}\ldots \e_{r_n} [s]$, we can easily prove that
$\Lambda(\alpha)$, when applied to the singleton $\{e\}$, gives
$\{r_1,\ldots,r_n,s,e\}$.

  Based on the existence of these two representations we may prove
uniqueness of our decomposition.

  \state Proposition \sl Every $\alpha$ in $\S G$ admits a unique
standard decomposition
    $$
  \alpha = \e_{r_1}\ldots \e_{r_n} [s]
    $$
  up to the order of the $\e_r$'s.

  \proof As observed above,
  $\Lambda (\alpha)(\{e\}) = \{r_1,\ldots,r_n,s,e\}$ and also
  $\degree(\alpha) = s$.
  So, if one has another standard decomposition
  $\alpha = \e_{t_1}\ldots \e_{t_m} [u]$ then, on the one hand, we would
have
  $s = u$
  and, on the other,
  $\{r_1,\ldots,r_n,s,e\} = \{t_1,\ldots, t_m,u,e\}$.
  Hence
    $$
  \{r_1,\ldots,r_n,s,e\} \backslash \{s,e\} =
  \{t_1,\ldots,t_m,u,e\} \backslash \{u,e\}
    $$
  which implies that
    $$
  \{r_1,\ldots,r_n\} =
  \{t_1,\ldots,t_m\} \proofend
    $$

  With this we can compute the order of $\S G$ for finite groups.

  \state Theorem \sl If $G$ is a finite group of order $p$, then $\S{G}$
has $2^{p-2}(p+1)$ elements.

  \proof In fact, there are $2^{p-1}$ elements in $\S{G}$ of the form
$\e_{r_1}\ldots \e_{r_n} [e]$ and for $p-1$ possible choices of $s\neq
e$, there are $2^{p-2}$ elements of the form $\e_{r_1}\ldots \e_{r_n}
[s]$.  The total number is then
  $2^{p-1} + (p-1)2^{p-2} = 2^{p-2}(p+1)$.
  \proofend

  Therefore we see that the order of $\S G$ grows exponentially with the
order of $G$.  For a group of order 28, for example, the order of
$\S{G}$ is 1,946,157,056.
  As for less dramatic examples, groups of order 2 trough 10 cause the
order of $\S{G}$ to be 3, 8, 20, 48, 112, 256, 576, 1280 and 2816,
respectively.

  \state Theorem \sl For every group $G$, $\S G$ is an \isg.

  \proof Assume $\alpha\in\S G$ admits, in addition to $\alpha^*$,
another ``inverse'', that is an element $\beta$ in $\S G$ such that
$\alpha\beta\alpha = \alpha$ and $\beta\alpha\beta = \beta$.  Write, in
standard form,
    $$
  \alpha =
  \e_{r_1}\ldots \e_{r_n} [s]
    $$
  and
    $$
  \beta^* =
  \e_{t_1}\ldots \e_{t_m} [u].
    $$
  Therefore $\beta =
  [u\inv]\e_{t_1}\ldots \e_{t_m}$ and we have
  $s =
  \degree(\alpha) =
  \degree(\alpha\beta\alpha) =
  s u\inv s$
  from which it follows that $u=s$. We also have
    $$
  \alpha\beta\alpha =
  \e_{r_1}\ldots \e_{r_n} [s]
  [s\inv]\e_{t_1}\ldots \e_{t_m}
  \e_{r_1}\ldots \e_{r_n} [s]
    \$ =
  \e_{r_1}\ldots \e_{r_n}
  \e_{t_1}\ldots \e_{t_m}[s] =
  \alpha.
    $$
  By the uniqueness of standard decompositions, we see that
    $$
  \{r_1,\ldots, r_n\} \cup \{t_1,\ldots, t_m\} =
  \{r_1,\ldots, r_n\}
    $$
  so
    $$
 \{t_1,\ldots, t_m\} \subseteq \{r_1,\ldots, r_n\}.
    $$
  The same argument, applied to the identity
  $\beta^*\alpha^*\beta^*=\beta^*$, yields
    $$
  \{r_1,\ldots, r_n\} \subseteq \{t_1,\ldots, t_m\},
    $$
  so we conclude that $\beta=\alpha^*$. \proofend

  \section {Actions of \isgs vs.~partial actions of groups} Recall that
a partial action of a group $G$ on a set $X$ is a map $\t \colon G \arw
\I(X)$ satisfying the conditions of \lcite{\PartialAction}.  Our next
proposition is designed to rephrase these conditions in terms of the \sg
structure of $\I(X)$.

  \nstate Proposition \ChangeLingua \sl Let $G$ be a group and $X$ be a
set. A map
  $\t \colon G \arw \I(X)$
  gives a partial action of $G$ on $X$ if and only if, for all $s,t\in
G$, we have
  \zitemno = 0
  \zitem $\t_s\t_t\t_{t\inv} =\t_{s t}\t_{t\inv}$,
  \zitem $\t_e = id_X.$
  \smallskip \noindent
  In this case $\t$ also satisfies
  \zitem $\t_{s\inv}\t_s\t_t =\t_{s\inv}\t_{s t}.$

  \proof Taking $s=t\inv$ in (i), we have
    $$
  \t_{t\inv}\t_t\t_{t\inv} =
  \t_e\t_{t\inv} =
  \t_{t\inv}.
    $$
  Replacing the roles of $t$ and $t\inv$, we also have
    $
  \t_t\t_{t\inv}\t_t =
  \t_t.
    $
  So, by the uniqueness of inverses in \isgs, we conclude that $\t_t^* =
\t_{t\inv}$.  Define $D_t = \ran(\t_t)$, so the above conclusion tells
us that
    $$
  \dom(\t_t) = \ran(\t_t^*) = \ran(\t_{t\inv}) = D_{t\inv}.
    $$
  That is, $\t_t$ is a map from $D_{t\inv}$ to $D_t$, as required.  Now,
for any $s,t\in G$ we have
    $$
  \t_{t\inv} \t_{s\inv} =
  \t_{t\inv} \t_{s\inv} \t_{s} \t_{s\inv} =
  \t_{t\inv s\inv} \t_{s} \t_{s\inv}.
    $$
  So, in particular, the domains of these maps must agree.  On the one
hand we have
    $$
  \dom(\t_{t\inv} \t_{s\inv}) = \t_s(D_{s\inv} \cap D_t).
    $$
  Noting that $\t_{s}\t_{s\inv}$ is the identity function on $D_s$, we
have, on the other hand,
    $$
  \dom(\t_{t\inv s\inv} \t_{s} \t_{s\inv}) = D_s \cap D_{s t}.
    $$
  Now the third condition of \lcite{\PartialAction} can easily be
checked.  The converse is left to the reader.

  With respect to the last part of the statement note that, since
  $\t_{t\inv}\t_{s\inv}\t_s = \t_{t\inv s\inv}\t_s$
  and
  $\t_t^*=\t_{t\inv}$, we have
    $$
  \t_{s\inv}\t_s\t_t =
  (\t_{t\inv}\t_{s\inv}\t_s)^* =
  (\t_{t\inv s\inv}\t_s)^* =
  \t_{s\inv}\t_{s t}. \proofend
    $$

  Our main result may now be stated.

  \state Theorem \sl For every group $G$ and any set $X$, there is a
one-to-one correspondence between
  \smallskip \item{(a)} partial actions of $G$ on $X$, and
  \smallskip \item{(b)} actions of $\S G$ on $X$.

  \proof By \lcite{\Universal}, \homos from $\S G$ into $\I(X)$ are in
one-to-one correspondence with maps from $G$ to $\I(X)$ satisfying
(i)--(iii) of \lcite{\Universal}.  As we have seen in
\lcite{\ChangeLingua}, these correspond to partial actions of $G$ on
$X$.
  \proofend

  \section {Graded $C^*$-algebras} References for $C^*$-algebra theory
are \cite{\Pedersen} and \cite{\FD}.

  \nstate Definition \Graded Let $G$ be a group.  A $C^*$-algebra $A$ is
said to be \stress{graded} over $G$, or $G$-\stress{graded}, if $A$ is
equipped with a linearly independent family of closed linear subspaces
$\{A_t\}_{t\in G}$, such that for all $s,t\in G$
  \zitemno = 0
  \zitem $A_s A_t\subseteq A_{st}$,
  \zitem $A_t^* = A_{t\inv}$,
  \zitem $A$ is the closure of $\directsum_{t\in G} A_t$.
  \smallskip \noindent In this case, each $A_t$ is called a
\stress{grading subspace} of $A$.

  We would like to show that $\S G$ plays an important role in
describing certain subspaces of $A$, related to its graduation.  For
example, given $s$ and $t$ in $G$, the product $A_s A_t$, although
contained in $A_{st}$, need not coincide with it, and, in fact, need not
even be dense there.  We could, therefore, ask what is the maximum
number of subspaces of $A$ one can obtain by taking the closure of
finite products of grading subspaces.

Before we proceed, let us agree on the following slightly unusual
notational convention: If $X$ and $Y$ are subsets of a $C^*$-algebra
then $XY$ will denote the \stress{closed} linear span of the set of
products $xy$ with $x \in X$ and $y \in Y$.

Returning to our graded algebra, let us denote by $D_t = A_tA_{t\inv}$.
It is easy to show that $A_t$ is a $D_t$-$D_{t\inv}$-Hilbert bimodule. A
good reference for the theory of Hilbert bimodules is \cite{\Jensen}.
As is the case for all right Hilbert modules, the product $A_t
D_{t\inv}$ must coincide with $A_t$ \scite{\Jensen}{1.1.4}.  This
translates to $A_tA_{t\inv}A_t = A_t$.

  \state Definition Given a $C^*$-algebra $A$, let $\Pl$ be the set of
all closed linear subspaces of $A$, equipped with the \sg structure
under which the product of $X$ and $Y$ is defined to be the closed
linear span of the set of products $xy$ with $x \in X$ and $y \in Y$, as
mentioned earlier.

  At this point one cannot oversee the fact that the map
    $$
  t\in G \mapsto A_t \in \Pl
    $$
  provides a key connection between the theory of graded $C^*$-algebras
and the theory of \sg\null s we have been discussing so far.  As already
mentioned, this map is not a \homo of \sg\null s but:

  \nstate Proposition \GradingProperties \sl Given a graded
$C^*$-algebra $A = \overline{\directsum_{t\in G} A_t}$, then, for every
$s$ and $t$ in $G$ one has
  \zitemno = 0
  \zitem $A_{s\inv} A_s A_t = A_{s\inv} A_{st}$,
  \zitem $A_s A_t A_{t\inv} = A_{st} A_{t\inv}$,
  \zitem $A_s A_e = A_s$.
  \smallskip

  \proof Based on the fact that $A_tA_{t\inv}A_t = A_t$ and that $A_s
A_t\subseteq A_{st}$ we have
    $$
  A_{s\inv} A_s A_t \subseteq
  A_{s\inv} A_{st} =
  A_{s\inv} A_s A_{s\inv} A_{st} \subseteq
  A_{s\inv} A_s A_t,
    $$
  which proves (i).  The proof of (ii) is quite similar, while (iii) can
be proved as follows:
    $$
  A_s A_e \subseteq
  A_s =
  A_s A_{s\inv} A_s \subseteq
  A_s A_e. \proofend
    $$

  \nstate Theorem \PartialGrading \sl Given a graded $C^*$-algebra $A =
\overline{\directsum_{t\in G} A_t}$, there is a correspondence which
assigns for each $\alpha$ in $\S G$, a closed subspace $A^\alpha$ of $A$
satisfying, for all $\alpha$ and $\beta$ in $\S G$ and all $t$ in $G$
  \zitemno=0
  \zitem $A^{[t]} = A_t$,
  \zitem if $\degree(\alpha) = t$ then $A^\alpha$ is contained in $A_t$,
  \zitem the closed linear span of the product of $A^\alpha$ by
$A^\beta$ is \stress{exactly} equal to $A^{\alpha\beta}$.
  \smallskip

  \proof This is an immediate consequence of \lcite{\GradingProperties}
and \lcite{\Universal}.
  \proofend

  As a consequence, we see that the collection $\{A^\alpha\}_{\alpha\in
\S G}$ is closed under multiplication, i.e, it is a sub\sg of $\Pl$.
Since it contains the $A_t$'s, we have provided an answer for the
question posed after \lcite{\Graded}, concerning the maximum number of
different subspaces of $A$ that one can obtain by taking finite products
of grading subspaces.  In the specific case of finite groups, that
number is at most the order of $\S G$, which we have seen to equal
$2^{p-2}(p+1)$, where $p$ is the order of $G$.

The question begging to be asked at this point is whether one can
produce a $G$-graded $C^*$-algebra such that the $A^\alpha$ are all
distinct.  In fact this question makes sense even for infinite groups.
The answer is yes but we are not yet prepared to provide an example
right now.  First we need to develop our theory a little further.

\section {Hilbert space representations} Contrary to our previous
treatment of representations of \sg\null s, simply as homomorphisms into
another \sg, we shall now make a more systematic study of
representations of \isgs on Hilbert spaces.  Let, therefore, $H$ denote
Hilbert's space and denote by $\B$ the collection of all bounded linear
operators on $H$.

\nstate Definition \HSRep Let $S$ be an \isg.  A \stress{(Hilbert space)
representation} of $S$ on $H$ is a map
  $\pi \colon S \arw \B$
  such that, for all $\alpha,\beta$ in $S$
  \zitemno =0
  \zitem $\pi(\alpha\beta) = \pi(\alpha)\pi(\beta)$,
  \zitem $\pi(\alpha^*) = \pi(\alpha)^*$.
  \smallskip\noindent If $S$ has a unit $e$, we shall also require that
  \zitem $\pi(e) = I$.

Note that, since $\alpha\alpha^*\alpha = \alpha$, the same holds for
$\pi(\alpha)$, which is then a partially isometric operator on $H$ for
each $\alpha \in S$.  Also, if $\e$ is an idempotent in $S$, one
necessarily has $\e^* = \e$ and hence $\pi(\e)$ is a selfadjoint
projection on $H$.

  The related notion for groups is that of a partial representation:

\state Definition Let $G$ be a group.  A \stress{partial representation}
of $G$ on $H$ is a map
  $\pi\colon G\arw \B$
  such that, for all $t,s\in G$,
  \zitemno = 0
  \zitem $\pi(s)\pi(t)\pi(t\inv) = \pi(st)\pi(t\inv)$,
  \zitem $\pi(t\inv) = \pi(t)^*$,
  \zitem $\pi(e) = I$.

  \smallskip Note that (i)--(iii) imply that $\pi(s\inv)\pi(s)\pi(t) =
\pi(s\inv)\pi(st)$.  So, regarding $\B$ as a \sg under multiplication
(i.e, composition of operators), we see that a partial representation of
$G$ satisfies the conditions of \lcite{\Universal} and hence there
exists a \homo\
    $$
  \tilde\pi \colon \S G \arw \B
    $$
  such that $\tilde\pi([t]) = \pi(t)$ for all $t$ in $G$.  It is not
hard to see that $\tilde\pi$ also satisfies $\tilde\pi(\alpha^*) =
\pi(\alpha)^*$, and hence $\tilde\pi$ is a representation of $\S G$ on
$H$ in the sense of \lcite{\HSRep}.

  Conversely, given a representation $\rho$ of $\S G$ on $H$, it is easy
to see that the map
    $$
  t \in G \mapsto \rho([t]) \in \B
    $$
  is a partial representation of $G$. We have therefore proven:

  \nstate Proposition \PreRep There is a one-to-one correspondence
between
  \smallskip \item{(a)} partial representations of $G$ on $H$, and
  \smallskip \item{(b)} representations of $\S G$ on $H$.

  \smallskip One of the important features of the theory of
$C^*$-algebras is that several different representation theories are
governed by a $C^*$-algebra.  The theory of partial representations of
groups is no exception and we shall now construct a $C^*$-algebra,
called the \stress{partial group $C^*$-algebra} of $G$ which will play
the corresponding role with respect to partial group representations.

The representation theory of \isgs has received much attention and
there, also, one can associate Banach and \cstar-algebras of interest to
the theory.  A key result, due to Barnes \cite{\Barnes}, states that the
$l_1$ algebra of an \isg has a separating family of irreducible
representations.  Also, Duncan and Paterson \cite{\DP} have studied both
the reduced and the full \cstar-algebras of \isg.

However, given our interest in the special case of the \isg associated
to a group, with all of its characteristic features, we shall prefer to
develop our theory from scratch.

  Let $G$ denote a group, which we will consider fixed for the remainder
of this section.

  An auxiliary \cstar-algebra in our study is the universal
\cstar-algebra, denoted $A$, defined via generators and relations as
follows.  The set of generators consists of a symbol $\p_E$ for each
finite subset $E\subseteq G$, and the relations are:
  \zitemno = 0
  \zitem $\p_E=\p_E^*$,
  \zitem $\p_E \p_F = \p_{E\cup F}$,
  \smallskip
  \noindent for all possible choices of $E$ and $F$.
  Note that each $\p_E$ will then be a projection and, if $\emptyset$
denotes the empty set, $\p_\emptyset$ will serve as a unit for $A$.
Obviously, $A$ is an abelian \cstar-algebra, so by Gelfand's Theorem,
$A=C(X)$ for some compact space $X$.  For the time being the description
of $A$ as a universal \cstar-algebra will be quite appropriate for our
purposes, so we shall postpone the study of its spectrum until needed.

  Since the collection of all $\p_E$'s is closed under multiplication
and adjoints, one has that $A$ is the closed linear span of these
elements.

  Let $A_e$ be the ideal of $A$ generated by $\p_{\{e\}}$, that is, $A_e
= \p_{\{e\}}A$. As before, $A_e$ is the closed linear span of the set
$\{\p_E\colon e\in E\}$. Also $A_e$ is a unital \cstar-algebra with unit
$\p_{\{e\}}$.

  For each $t$ in $G$, let $\t_t$ be the automorphism of $A$ specified
by $\t_t(\p_E) = \p_{tE}$.  Observe that $A_e$ is not invariant under
$\t_t$, but if we let
    $$
  D_t = \span\{\p_E\colon e,t\in E\},
    $$
  then $D_t$ is an ideal in $A_e$ and the restriction of $\t_t$ to
$D_{t\inv}$ is a \cstar-algebra isomorphism
    $$
  \t_t\colon D_{t\inv} \arw D_t
    $$
  which lives fully inside of $A_e$. This is an example of a partial
automorphism of $A_e$, as defined in \cite{\ExelCircle}, and the map
  $t \mapsto \t_t$
  fits the definition of a partial action of $G$ on $A_e$ given in
\cite{\ExelTPA} and \cite{\McClanahan}.
  In particular, we can form the crossed product of $A_e$ by $G$.

  \state Definition The \stress{partial group $C^*$-algebra} of $G$ is
the $C^*$-algebra $C^*_p(G)$ given by the crossed product
\cite{\ExelTPA, \McClanahan} of $A_e$ by $\t$, that is,
    $$
  C^*_p(G) = A_e\crossproduct_\t G.
    $$

  The following result improves upon \lcite{\PreRep} and shows the
relevance of $C^*_p(G)$ in the theory of partial representations of $G$.

  \state Theorem \sl There is a one-to-one correspondence between
  \smallskip \item{(a)} partial representations of $G$ on $H$,
  \smallskip \item{(b)} representations of $\S G$ on $H$ and
  \smallskip \item{(c)} \cstar-algebra representations of $C^*_p(G)$ on
$H$.

  \proof Let $r,s$ be in $G$ and let $E$ and $F$ be finite subsets of
$G$ containing $\{e,r\}$ and $\{e,s\}$, respectively.  Then $\p_E$ is in
$D_r$ and $\p_F$ is in $D_s$ and we have, following \cite{\ExelTPA} and
\cite{\McClanahan},
    $$
  (\p_E\delta_r) * (\p_F\delta_s) =
  \t_r ( \t_r\inv(\p_E)\p_F )\delta_{rs} =
  \t_r ( \p_{r\inv E\cup F} )\delta_{rs} =
  \p_{ E\cup r F} \delta_{rs}.
    $$

Consider, for each $t$ in $G$, the element $u_t = \p_{\{e,t\}}\delta_t
\in C^*_p(G)$.  Based on the computation above, it is easy to verify
that $u_s u_t u_{t\inv} = u_{st}u_{t\inv}$, that $u_{t\inv} = u_t^*$ and
that $u_e = 1$.  This said, it is easy to see that every representation
of $C^*_p(G)$, when computed on the elements $u_t$, gives a partial
representation of $G$.

  Conversely, suppose we are given a representation $\pi$ of $\S G$ on
$H$.  Recall that the idempotents $\e_r$ in $\S G$ are defined by $\e_r
= [r][r\inv]$.  For each finite set $E = \{r_1,\ldots,r_n\} \subseteq
G$, define
    $$
  \q_E = \pi(\e_{r_1}\ldots \e_{r_n}).
    $$
  This is a projection in $\B$ and we have $\q_E \q_F = \q_{E \cup F}$
for all $E$ and $F$, so there exists a representation $\rho$ of $A$ on
$H$ such that
  $\rho(\p_E) = \q_E$.  Also let $u_t = \pi([t])$.

  Observe that, for any $t,r$ in $G$ one has
    $$
  [t]e_r =
  [t][r][r\inv] =
  [tr][r\inv] =
  [tr] [r\inv t\inv] [tr] [r\inv]
  \$=
  [tr] [r\inv t\inv] [t] =
  \e_{tr} [t].
    $$

  \noindent Therefore, if $E$ contains $\{e,t\inv\}$, then
    $$
  u_t \rho(\p_E) u_{t\inv} =
  \pi([t] \e_{r_1}\ldots \e_{r_n} [t\inv] ) =
  \pi(\e_{tr_1}\ldots \e_{tr_n}[t] [t\inv] )
  \$=
  \pi(\e_{tr_1}\ldots \e_{tr_n} ) =
  \rho(\p_{tE}) =
  \rho(\t_t(\p_{E})).
    $$
  So we get a covariant representation of the partial action of $G$ on
$A_e$, in the sense of \scite{\McClanahan}{2.8}, and hence a
representation $\rho \x u$ of $C^*_p(G)$.  In particular
    $$
  \rho \x u (\p_{\{e,t\}}\delta_t) =
  \rho (\p_{\{e,t\}})u_t =
  \q_{\{e,t\}} \pi([t])=
  \pi(\e_e\e_t[t]) = \pi([t]). \proofend
    $$

  For our following result it will be important to understand the
algebra $A$ more concretely.  For this purpose let $K$ be the power set
of $G$.  Among many different descriptions for $K$, one can use the
model
    $$
  K = \prod_{t \in G} \{0,1\}.
    $$

  Equipped with the product topology, $K$ becomes a compact Hausdorff
space.  When $G$ is an infinite countable group, for example, $K$ is
homeomorphic to Cantor's set.  For each $t$ in $G$ let $\q_t\in C(K)$ be
the $t^{th}$ coordinate function, i.e,
    $$
  \q_t \colon x\in K \mapsto x_t \in \C.
    $$
  For each finite set $E\subseteq G$ let $\q_E = \prod_{t\in E} \q_t$.
Clearly $\q_E \q_F = \q_{E \cup F}$ so we get a \cstar-algebra
homomorphism
    $$
  \Phi \colon A \arw C(K)
    $$
  such that $\Phi(\p_E) = \q_E$.

  \state Proposition \sl $A$ and $C(K)$ are isomorphic via the
homomorphism $\Phi$ above.

  \proof Surjectivity is an immediate consequence of the
Stone--Weierstrass Theorem.  To prove injectivity let $a\in A$ be non-z
ero.  Take a character $f$ of $A$ such that $f(a)\neq 0$. For each $t$
in $G$ observe that $\p_{\{t\}}$ is a projection in $A$, so
$f(\p_{\{t\}})$ is a complex number in $\{0,1\}$.  Define $x$ in $K$ by
    $$
  x = \bigl( f(\p_{\{t\}}) \bigr)_{t\in G}.
    $$
  For all $t$ we have $\Phi(\p_{\{t\}})(x) = \q_t(x) = f(\p_{\{t\}})$
so, observing that the $\p_{\{t\}}$ generate $A$, we see that we must
also have $\Phi(a)(x) = f(a) \neq 0$, from where we conclude that
$\Phi(a) \neq 0$ and hence that $\Phi$ is injective. \proofend

  It is well known that a crossed product by a partial action is a
graded $C^*$-algebra, and hence so is $C^*_p(G)$.  This is also the case
for the full (as opposed to partial) group $C^*$-algebra $C^*(G)$ (see
\cite\Pedersen, for definitions).  However the grading of $C^*(G)$ is
rather uninteresting from our point of view since it is saturated, that
is, the product of the grading subspaces $A_t$ and $A_s$ always
coincides with $A_{ts}$.  In terms of the notation of
\lcite{\PartialGrading}, this simply says that $A^\alpha =
A_{\degree(\alpha)}$.  The situation for the partial group algebra is a
lot more interesting.

  \nstate Theorem \Biggrading \sl Let $C^*_p(G) =
\overline{\directsum_{t\in G} B_t}$ be the $G$-grading of $C^*_p(G)$
arising from its construction as a partial crossed product.  Also let
    $$
  \alpha \in \S G \mapsto B^\alpha \in \Pl
    $$
  be the \sg-\homo referred to in \lcite{\PartialGrading}.  Then, if
$\alpha \neq \beta$, one has $B^\alpha \neq B^\beta$.

  \proof We initially claim that if $\alpha = \e_{r_1}\ldots \e_{r_n}[s]
\in \S G$. Then
    $$
  B^\alpha = \span\{\p_E\delta_s\colon e,r_1,\ldots,r_n,s \in E\}.
    $$
  In fact, if $n=0$, then $B^\alpha = B^{[s]} = B_s$ and the result
follows from the definition of $B_s$ as $B_s = D_s \delta_s$ and the
fact that $D_s$ is spanned by the $\p_E$'s with $e,s\in E$.

  Let $n>0$ and $\gamma=\e_{r_2}\ldots \e_{r_n}[s]$ so
$\alpha=\e_{r_1}\gamma$ and $B^\alpha = B^{[r_1][r_1\inv]\gamma} =
B_{r_1}B_{r_1\inv}B^\gamma$.  By induction we may assume that $B^\gamma
= \span\{\p_E\delta_s\colon e,r_2,\ldots,r_n,s \in E\}$.

   On the other hand, it is easy to see that $B_{r_1}B_{r_1\inv} =
D_{r_1} = \span \{\p_F\colon e,r_1 \in F\}$.  We thus have
    $$
  B_\alpha = \span \{\p_F\p_E\delta_s\colon e,r_1 \in F ;
e,r_2,\ldots,r_n,s \in E\}.
    $$
  Observing that $\p_F\p_E = \p_{F\cup E}$ we see that the claim is
proved.

  Returning to the proof of our Theorem, suppose
    $$
  \alpha = \e_{r_1}\ldots \e_{r_n} [s]
    $$
  and
    $$
  \beta = \e_{t_1}\ldots \e_{t_m} [u].
    $$
  are in standard form and that $B^\alpha = B^\beta$.  Then, because of
\lcite{\PartialGrading.ii}, that is, because
  $B^\alpha \subseteq B_s$
  and
  $B^\beta\subseteq B_u$,
  we must have $s=u$.

 On the other hand, $B^\alpha B^{\alpha^*} = B^{\alpha\alpha^*}$, but
  $\alpha\alpha^* = \e_{r_1}\ldots \e_{r_n} \e_s$. So, by our claim, we
have
  $B^{\alpha\alpha^*} = \span\{\p_E\delta_e\colon e,r_1,\ldots,r_n,s \in
E\}$,
  which is a subspace of $B_e$. Since $B^\alpha = B^\beta$, this
subspace must coincide with
  $B^{\beta\beta^*} = \span\{\p_E\delta_e\colon e,t_1,\ldots,t_m,u \in
E\}$ from which the final conclusion follows easily.
  \proofend

  In the special case of finite groups we have therefore proven:

  \state Corollary \sl Let $G$ be a finite group of order $p$.  Then,
there is a $G$-graded $C^*$-algebra $A = \directsum_{t\in G} A_t$, such
that the sub\sg of $\Pl$ generated by the grading spaces has exactly
$2^{p-2}(p+1)$ elements; while, for any $G$-graded $C^*$-algebra the
said sub\sg cannot have more than this number of elements.

  Beyond the purpose of constructing interesting graded $C^*$-algebras,
we believe there is a considerable amount of interesting information
that can be dug out from the structure of $C^*_p(G)$.

  To deal with the specific case of finite groups, one can prove that
$C^*_p(G)$ is also the groupoid $C^*$-algebra \cite{\Renault}
constructed from the partial action of $G$ on $C(K)$, and compute
$C^*_p(G)$ for the most elementary finite groups. In particular, for the
two groups of order 4, namely $Z/4Z$ and $Z/2Z \+ Z/2Z$ one obtains
    $$
  C^*_p(Z/4Z) = C^7 \+ M_2(C) \+ M_3(C)
    $$
  while
    $$
  C^*_p(Z/2Z \oplus Z/2Z) = C^{11} \+ M_3(C).
    $$
  It is quite intriguing that, while the group $C^*$-algebra of finite
commutative groups forgets everything but the order of the group, the
partial group $C^*$-algebra $C^*_p(G)$ allows us to distinguish between
$Z/4Z$ and $Z/2Z \+ Z/2Z$.

  \bigbreak
  \centerline{\tensc References}
  \nobreak\medskip\frenchspacing

\ATtechreport{\AEE,
  author = {B. Abadie, S. Eilers and R. Exel},
  title = {Morita Equivalence and Crossed Products by {Hilbert}
Bimodules},
  institution = {Universidade de S\~ao Paulo},
  year = {1994},
  note = {preprint},
  toappear = {Trans. Amer. Math. Soc},
  NULL = {},
  atrib = {IR},
  MR = {}
  }

\ATarticle{\Barnes,
  author = {B. A. Barnes},
  title = {Representations of the $l_1$ algebra of an inverse
semigroup},
  journal = {Trans. Amer. Math. Soc},
  year = {1976},
  volume = {218},
  pages = {361--396},
  NULL = {},
  }

\ATarticle{\DP,
  author = {J. Duncan and L. T.  Paterson},
  title = {$C^*$-algebras of inverse semigroups},
  journal = {Proc. Edinburgh Math. Soc.},
  year = {1985},
  volume = {28},
  pages = {41-58},
  NULL = {},
  }

\ATtechreport{\ExelTPA,
  author = {R. Exel},
  title = {Twisted Partial Actions, A Classification of Stable
{$C^*$}-Algebraic Bundles},
  institution = {Universidade de S\~ao Paulo},
  year = {1994},
  note = {preprint},
  toappear = {},
  NULL = {},
  atrib = {N}
  }

\ATarticle{\ExelCircle,
  author = {R. Exel},
  title = {Circle Actions on {$C^*$}-Algebras, Partial Automorphisms and
a Generalized {Pimsner}--{Voiculescu} Exact Sequence},
  journal = {J. Funct. Analysis},
  year = {1994},
  volume = {122},
  pages = {361--401},
  NULL = {},
  atrib = {IR},
  MR = {95g:46122}
  }

\ATbook{\FD,
  author = {J. M. G. Fell and R. S. Doran},
  title = {Representations of *-algebras, locally compact groups, and
Banach *-algebraic bundles},
  publisher = {Academic Press},
  year = {1988},
  volume = {125 and 126},
  series = {Pure and Applied Mathematics},
  NULL = {},
  }

\ATbook{\Howie,
  author = {J. M. Howie},
  title = {An introduction to Semigroup theory},
  publisher = {Academic Press},
  year = {1976},
  volume = {},
  series = {},
  NULL = {},
  }

\ATbook{\Jensen,
  author = {K. Jensen and K. Thomsen},
  title = {Elements of $K\!K$-Theory},
  publisher = {Birkh\"auser},
  year = {1991},
  volume = {},
  series = {},
  NULL = {},
  }

\ATtechreport{\McClanahan,
  author = {K. McClanahan},
  title = {$K$-theory for partial crossed products by discrete groups},
  institution = {University of Mississippi},
  year = {1994},
  note = {preprint},
  toappear = {},
  NULL = {},
  }

\ATbook{\Pedersen,
  author = {G. K. Pedersen},
  title = {$C^*$-Algebras and their automorphism groups},
  publisher = {Acad. Press},
  year = {1979},
  volume = {},
  series = {},
  NULL = {},
  }

\ATtechreport{\Pimsner,
  author = {M. Pimsner},
  title = {A class of $C^*$-algebras generalizing both Cuntz-Krieger
algebras and crossed products by {\bf Z}},
  institution = {University of Pennsylvania},
  year = {1993},
  note = {preprint},
  toappear = {},
  NULL = {},
  }

\ATtechreport{\QuiggRaeburn,
  author = {J. Quigg and I. Raeburn},
  title = {Landstad duality for partial actions},
  institution = {University of Newcastle},
  year = {1994},
  note = {preprint},
  toappear = {},
  NULL = {},
  }

\ATbook{\Renault,
  author = {J. Renault},
  title = {A groupoid approach to $C^*$-lgebras},
  publisher = {Springer},
  year = {1980},
  volume = {793},
  series = {Lecture Notes in Mathematics},
  NULL = {},
  }

\ATmasterthesis{\Sieben,
  author = {N. Sieben},
  title = {$C^*$-crossed products by partial actions and actions of
inverse semigroups},
  school = {Arizona State University},
  year = {1994},
  NULL = {},
  }

  \vskip 2cm
  \rightline{November 1995}

  \bye